\documentclass[12pt]{amsart}
\usepackage{graphicx}
\newtheorem{theorem}{Theorem}[section]

\newtheorem{lemma}[theorem]{Lemma}

\newtheorem{corollary}[theorem]{Corollary}

\theoremstyle{remark}

\newtheorem{conjecture}[theorem]{Conjecture}

\newtheorem{definition}[theorem]{Definition}

\numberwithin{equation}{section}

\begin{document}

\title[
Canonical forms for operation tables of finiate connected quandles 
]{
Canonical forms for operation tables of finite connected quandles
}

\author{Chuichiro Hayashi}

\date{\today}

\thanks{The first author is partially supported
by Grant-in-Aid for Scientific Research (No. 22540101),
Ministry of Education, Science, Sports and Technology, Japan.}

\begin{abstract}
 We introduce a notion of natural orderings of elements of finite connected quandles of order $n$.
 When the elements of such a quandle $Q$ are already ordered naturally,
any automophism on $Q$ is a natural ordering.
 Although there are many natural orderings,
the operation tables for such orderings coincide
when the permutation $*q$ is a cycle of length $n-1$.
 This leads to the classification of automorphisms on such a quandle.
 Moreover, it is also shown 
that every row and column of the operation table of such a quandle contains all the elements of $Q$,
which is due to K. Oshiro.
 We also consider the general case of finite connected quandles.
\end{abstract}

\maketitle

\section{Introduction}

 The algebraic structure of quandle was introduced
by D. Joyce and S. V. Matveev in \cite{J} and \cite{M}.
 There was defined an invariant for a classical knot called knot quandle,
which classifies knots up to homeomorphism of pairs.

 A {\it quandle} is a set $Q$ with a binary operation $*$ : $Q \times Q \rightarrow Q$
satisfying the three axioms
\begin{enumerate}
\item[(1)]
for any $a \in Q$, $a*a=a$,
\item[(2)]
for any pair $a,b \in Q$, 
there exists a unique $c \in Q$ such that $a = c*b$, and 
\item[(3)]
for any triple $a,b,c \in Q$, we have $(a*b)*c = (a*c)*(a*c)$.
\end{enumerate}

 Note that possibly $a*b \ne b*a$ and $a*(b*c) \ne (a*b)*(a*c)$ for some $a,b,c \in Q$.
 Axiom (3) is called right-distributivity.
 Axiom (2) is called right-invertibility,
and implies
that the map $r_b : Q \ni x \mapsto x*b \in Q$ is a bijection for all $b \in Q$.
 A quandle is called {\it trivial}
if $r_b$ is the identity map for all $b \in Q$.

 There is an inverse map $r_b^{-1}$,
and we denote $r_b^{-1} (a)$ by $a \bar{*} b$.
 Then $\bar{*}$ gives a binary operation on $Q$,
and under this operation $Q$ forms a quandle $(Q,\bar{*})$,
which we call the {\it dual} quandle of $(Q,*)$.
 The three formulae below are well-known.

\begin{enumerate}
\item[(4)]
For any triple $a,b,c \in Q$, $(a*b)\bar{*}c=(a\bar{*}c)*(b\bar{*}c)$. 
\item[(5)]
For any triple $a,b,c \in Q$, $(a\bar{*}b)*c=(a*c)\bar{*}(b*c)$. 
\item[(6)] 
For any triple $a,b,c \in Q$, $a*(b*c) = (((a\bar{*}c)*b)*c)$.
\end{enumerate}

 For example, any group $G$ forms a quandle under conjugation,
i.e., the operation $*$ defined by $a*b = b^{-1} a b$.
 Such a quandle is called the {\it conjugation quandle} of $G$.

 Let $Q_1, Q_2$ be quandles.
 A map $f : Q_1 \rightarrow Q_2$ is said to be a {\it homomorphism}
if $f(a*b)=f(a)*'f(b)$ holds for all $a,b \in Q_1$,
where $*$ and $*'$ are quandle operations in $Q_1$ and $Q_2$ respectively.
 If such a map $f$ is bijective,
then it is called an {\it isomorphism},
and we say that $Q_1$ and $Q_2$ are {\it isomorphic}. 
 An isomorphism from a quandle $Q$ to $Q$ itself
is called an {\it automorphism} of $Q$.
 Axiom (3) implies 
that the above bijection $r_b : Q \ni x \mapsto x*b \in Q$
is an automorphism of $Q$ for any $b \in Q$.
 Actually, we can rewrite (3)
as $r_c(a*b) = (a*b)*c = (a*c)*(a*c) = r_c(a) * r_c(b)$.

 By distinguishing quandles, we can distinguish knots.
 However, knot quandles have infinite number of elements,
and it is a hard problem
to decide given two knot quandles are isomorphic or not.
 Considering homomorphism from a knot quandle to a finite quandle
gives a convenient way to distinguish knot quandles.
 See, for example, \cite{FR}, \cite{CJKLS}, \cite{L} and \cite{DL}.
 For example,
the involutory knot quandle
is the knot quandle $Q$ 
with the condition $(a*b)*b=a$ for any pair $a, b \in Q$ added,
and is finite for many knots.
 See Section 19 in \cite{J}.

 A quandle $Q$ is said to be connected (or indecomposable)
if the orbit
\newline
$O(a) = \{ (\cdots ((a \star_1 b_1) \star_2 b_2) \star_3 \cdots \star_m b_m) 
\ |\ b_i \in Q,\ \star_i \in \{ *, \bar{*} \},\ m \in \{ 0,1,2,\cdots \} \}$
\newline
is equal to $Q$ for all $a \in Q$.
 Knot quandles are known to be connected.
 L. Vendramin classified connected quandles with $35$ or smaller number of elements in \cite{V}.

\begin{table}
\begin{center}
\begin{tabular}{c|cccccc}
  & 1 & 2 & 3 & 4 & 5 & 6 \\
\hline
1 & 1 & 1 & 5 & 6 & 3 & 4 \\
2 & 2 & 2 & 6 & 5 & 4 & 3 \\
3 & 5 & 6 & 3 & 3 & 1 & 2 \\
4 & 6 & 5 & 4 & 4 & 2 & 1 \\
5 & 3 & 4 & 1 & 2 & 5 & 5 \\
6 & 4 & 3 & 2 & 1 & 6 & 6 \\
\end{tabular}
\end{center}
\caption{\rm An operation table for the quandle $Q_{61}$}
\label{table:q61}
\end{table}

 Let $Q$ be a finite quandle of order $n$,
with its $n$ elements ordered, say, $q_1, q_2,\cdots, q_n$.
 We write them $1,2,\cdots, n$ for short.
 The {\it quandle matrix} for $Q$ is an $n \times n$ matrix
whose element $q_{ij}$ in row $i$ and column $j$ is $i*j$.
 Note that the diagonal element $q_{ii}$ is equal to $i$ by Axiom (1).
 The operation table for the quandle $Q_{61}$ in Vendramin's list is displayed in Table \ref{table:q61}.
 In the right bottom $6 \times 6$ quandle matrix for $Q_{61}$, 
we can see $3*4 = 3$ from the element in row $3$ and column $4$, for example.
 There are some papers on quandle matrices.
 See, for example, \cite{HN}, \cite{NW} and \cite{MNT}.
 Let $\nu : \{ 1,2, \cdots, n \} \rightarrow \{ 1,2, \cdots, n \}$ be a bijection.
 After reordering the elements of $Q$ by $\nu$,
the new quandle matrix
has the element $\nu(i*j)$ 
in row $\nu (i)$ and column $\nu (j)$ (see Definition 1 in \cite{HN}),
and $\nu$ is an automorphism of $Q$
if and only if the quandle matrix is unchanged under $\nu$ (Corollary 5 in \cite{HN}).
 For example, the automorphism $r_b$ does not change the quandle matrix
for any $b \in Q$.

 Since the map $r_j : Q \ni x \mapsto x*j \in Q$ is a bijection,
it can be regarded as a permutation on the set $\{ 1,2, \cdots, n \}$.
 We can see the permutation $r_j$ 
in the $j$th column of the quandle matrix for $Q$.
 For instance,
the $4$th column in the quandle matrix in Table \ref{table:q61}
shows that the operation $*4$ gives a permutation
$r_4 = 
\left( \begin{array}{ccccccc}
1 & 2 & 3 & 4 & 5 & 6 \\
6 & 5 & 3 & 4 & 2 & 1 \\
\end{array} \right)
=
(1\ 6)(2\ 5)(3)(4)$,
where
$(1\ 6)$, $(2\ 5)$, $(3)$ and $(4)$ are mutually disjoint cycles.

 As is well-known,
any permutation $\sigma$ is decomposed into a product of disjoint cycles
uniquely modulo ordering of the cycles.
 When
$\sigma = 
(i_{1,1}\ \cdots\ i_{1,\ell_1})
(i_{2,1}\ \cdots\ i_{2,\ell_2})
\cdots
(i_{k,1}\ \cdots\ i_{k,\ell_k})$,
we call the multiple set of the length of the cycles $\{ \ell_1, \ell_2, \cdots, \ell_k \}$ 
the {\it pattern} of $\sigma$,
where a multiple set admits repeats
and disregards ordering of its elements.
 For example, the pattern of the above permutation $r_4$ is $\{ 1,1,2,2 \}$.
 Note that 
the patterns of two permutations $\sigma$ and $\rho$ coincide
if and only if $\sigma$ and $\rho$ are mutually conjugate,
i.e., there is some permutation $\omega$ with $\rho = \omega^{-1} \sigma \omega$.

 In \cite{LR},
P. Lopes and D. Roseman defined the {\it profile} of a quandle with $n$ elements
to be the sequence of the patterns of $r_1, r_2, \cdots, r_n$.
 In case of a connected quandle $Q$,
$r_i$ and $r_j$ are mutually conjugate for any pair $i, j$ with $1 \le i < j \le n$
(Corollary 5.1 in \cite{LR}).
 This can be easily seen using Formula (6).
 Hence we call the pattern of $r_n$ the profile of $Q$ for short in this paper.
 For example, the profile of $Q_{61}$ given in Table \ref{table:q61} is $\{ 1,1,2,2 \}$.
 In general, a non-trivial finite connected quandle is of profile $\{ 1, \ell_1, \ell_2, \cdots, \ell_k \}$
with $1 \le \ell_1 \le \ell_2 \le \cdots \le \ell_k$ for some $k \in {\mathbb N}$.
We have at least one $1$ in profile since $r_i (i) = i$.

\begin{conjecture}
 $\ell_k$ is a multiple of $\ell_i$ for any integer $i$ with $1 \le i \le k-1$.
\end{conjecture}

 Lopes and Roseman studied
finite quandles with profile $(\{ 1, n-1 \}, \{ 1, n-1 \}, \cdots, \{ 1, n-1 \})$
in Theorem 6.5 and Corollaries 6.4--6.8 in \cite{LR}.
 They showed
that the $n$th permutation $r_n$ is $(1\ 2\ \cdots n-1)(n)$ modulo isomorphism, 
$r_{n-1}$ is a solution to a certain system of equations,
and $r_k$ with $1 \le k \le n-2$ is determined
by the formula $r_k = r_n^k r_{n-1} r_n^k$.

\begin{table}
\begin{center}
\begin{tabular}{c|ccccccc}
  & 1 & 2 & 3 & 4 & 5 & 6 & 7 \\
\hline
1 & 1 & 5 & 2 & 6 & 3 & 7 & 4 \\
2 & 5 & 2 & 6 & 3 & 7 & 4 & 1 \\
3 & 2 & 6 & 3 & 7 & 4 & 1 & 5 \\
4 & 6 & 3 & 7 & 4 & 1 & 5 & 2 \\
5 & 3 & 7 & 4 & 1 & 5 & 2 & 6 \\
6 & 7 & 4 & 1 & 5 & 2 & 6 & 3 \\
7 & 4 & 1 & 5 & 2 & 6 & 3 & 7 \\
\end{tabular}
\end{center}
\caption{A quandle matrix for $Q_{72}$}
\label{table:q72}
\end{table}

 We define natural reorderings of the elements of a finite connected quandle.

\begin{definition}\label{def:natural}
 Let $Q$ be a connected quandle with $n$ elements $1,2,\cdots,n$.
 A bijection $\nu : \{ 1,2, \cdots, n \} \rightarrow \{ 1,2,\cdots, n \}$
is called a {\it natural reordering} (with respect to $r_q$)
if $\nu (i_{st})= (\sum_{j=1}^{s-1} \ell_j) + t$ and $\nu (q)=n$
for some element $q \in Q$ and some presentation of the permutation $r_q$
as a product of disjoint cycles
$r_q = 
( i_{11}\ i_{12}\ \cdots \ i_{1\ell_1} )
( i_{21}\ i_{22}\ \cdots \ i_{2\ell_2} )
\cdots 
( i_{k1}\ i_{k2}\ \cdots \ i_{k\ell_k} )
(q)$
with $1 \le \ell_1 \le \ell_2 \le \cdots \le \ell_k$.
\end{definition}

 For example,
for the quandle $Q_{72}$ shown in Table \ref{table:q72},
$r_1$ is decomposed to $(2\ 5\ 3)(4\ 6\ 7)(1)$. 
 Hence
$\lambda(1)=7, \lambda(2)=1, \lambda(3)=3, \lambda(4)=4, \lambda(5)=2, \lambda(6)=5$ and $\lambda(7)=6$
give a reordering 
$\lambda : \{ 1,2,\cdots,7 \} \rightarrow \{ 1,2,\cdots,7\}$
which is natural with respect to $r_1$.
 Since we can rewrite $r_1$ as $(4\ 6\ 7)(2\ 5\ 3)(1)$,
$\mu(1)=7, \mu(2)=4, \mu(3)=6, \mu(4)=1, \mu(5)=5, \mu(6)=2$
and $\mu(7)=3$ also give a natural reordering $\mu$.
 The presentation $r_1 = (6\ 7\ 4)(3\ 2\ 5)(1)$ gives
another natural reordering $\nu$ with
$\nu(1)=7, \nu(2)=5, \nu(3)=4, \nu(4)=3, \nu(5)=6, \nu(6)=1$ and $\nu(7)=2$.
 Since $r_2=(1\ 5\ 7)(3\ 6\ 4)(2)$, 
$\xi(1)=1, \xi(2)=7, \xi(3)=4, \xi(4)=6, \xi(5)=2, \xi(6)=5$ and $\xi(7)=3$
determine a natural reordering $\xi$ with respect to $r_2$.

 After the elements of $Q$ are reordered naturally,
the permutation $r_n$ is decomposed
into disjoint cycles in the form below,
where $1 \le \ell_1 \le \ell_2 \le \cdots \le \ell_k$.
 We call this form (N).
\newline 
$r_n=(1\ 2\ \cdots \ell_1)
(\ell_1 + 1\ \ \ell_1 + 2\ \ \cdots\ \ \ell_1 + \ell_2)
(\ell_1 + \ell_2 + 1\ \ \ell_1 + \ell_2 + 2\ \ \cdots\ \ \ell_1 + \ell_2 + \ell_3)$
\newline
\hspace*{10mm}
$\cdots
((\sum_{j=1}^{k-1} \ell_j)+1\ \ (\sum_{j=1}^{k-1} \ell_j)+2\ \ \cdots \ \ (\sum_{j=1}^{k-1} \ell_j)+\ell_k) 
(n)$
\newline
 Note that there are $n\ell_1 \ell_2 \cdots \ell_k$ natural reorderings
if $1 \le \ell_1 < \ell_2 < \cdots < \ell_k$.
 There are more if $\ell_i = \ell_{i+1}$ for some $1 \le i \le k-1$.

\begin{definition}
 Let $Q$ be a finite connected quandle with $n$ elements $1,2,\cdots, n$. 
 If $r_n$ is decomposed into the form (N),
then we say that the elements of $Q$ are {\it naturally ordered},
and that the quandle matrix of $Q$ is in {\it canonical form}.
\end{definition}

\begin{theorem}\label{thm:automorphism}
 Let  $Q$ be a finite connected quandle whose elements are naturally ordered.
 Then any automorphism of $Q$ is a natural reordering.
\end{theorem}

\begin{theorem}\label{thm:canonical}
 Let $Q$ be a finite connected quandle 
with its elements $1,2,\cdots,n$ naturally ordered,
and with profile $\{ 1, \ell_1, \ell_2, \cdots, \ell_k \}$,
where $1 \le \ell_1 \le \ell_2 \le \cdots \le \ell_k$.
 For any reordering $\mu$ which is natural with respect to $r_q$
for some $q \in Q$,
there is a natural reordering $\nu$ with respect to $r_n$
such that
the quandle matrices after reordering by $\mu$ and $\nu$ coincide.
 Moreover, we can take $\nu$
so that, for some integer $m$ with $1 \le m \le k$ and $\ell_m = \ell_k$,
$\nu((\sum_{j=1}^{m-1} \ell_j)+i)=(\sum_{j=1}^{k-1} \ell_j)+i$ 
for any integer $i$ with $1 \le j \le \ell_m$.
 In particular, when $\ell_{k-1} < \ell_k$, 
we can take $\nu$
so that $\nu(i)=i$ 
for any integer $i$ with $(\sum_{j=1}^{k-1} \ell_j) + 1 \le i \le n$.
\end{theorem}

\begin{conjecture}
 In the last sentence of Theorem \ref{thm:canonical}, 
we can take $\nu$ 
so that $\nu(i)=i$ 
for any integer $i$ with $(\sum_{j=1}^{k-1} \ell_j) + 1 \le i \le n$
even when $\ell_{k-1} = \ell_k$.
\end{conjecture}

\begin{corollary}\label{cor:AtMost}
 Let $Q$ be a finite connected quandle 
with its elements naturally ordered,
and with profile $\{ 1, \ell_1, \ell_2, \cdots, \ell_k \}$,
where $1 \le \ell_1 < \ell_2 < \cdots < \ell_k$.
 Then the number of canonical forms of quandle matrices for $Q$ after reorderings
are at most $\ell_1 \ell_2 \cdots \ell_{k-1}$.
 In particular, when $k=1$ or $\lq\lq k=2$ and $\ell_1=1 "$,
the canonical form of quandle matrix for $Q$ is unique.
\end{corollary}

 The canonical forms of the quandle matrices for $Q_{52}$ and $Q_{53}$ in Vendramin's list
are shown in Tables \ref{table:CanonicalQ52} and  \ref{table:CanonicalQ53}.
 They are of profile $\{ 1,4 \}$.
 We can see that $Q_{52}$ and $Q_{53}$ are not isomorphic
because the canonical forms of their quandle matrices are distinct.

\begin{table}
\begin{minipage} {0.48\hsize}
\begin{center}
\begin{tabular}{c|ccccc}
  & 1 & 2 & 3 & 4 & 5 \\
\hline
1 & 1 & 4 & 5 & 3 & 2 \\
2 & 4 & 2 & 1 & 5 & 3 \\
3 & 5 & 1 & 3 & 2 & 4 \\
4 & 3 & 5 & 2 & 4 & 1 \\
5 & 2 & 3 & 4 & 1 & 5 \\
\end{tabular}
\end{center}
\caption{The canonical quandle matrix for $Q_{52}$}
\label{table:CanonicalQ52}
\end{minipage}
\begin{minipage} {0.48\hsize}
\begin{center}
\begin{tabular}{c|ccccc}
  & 1 & 2 & 3 & 4 & 5 \\
\hline
1 & 1 & 5 & 4 & 3 & 2 \\
2 & 4 & 2 & 5 & 1 & 3 \\
3 & 2 & 1 & 3 & 5 & 4 \\
4 & 5 & 3 & 2 & 4 & 1 \\
5 & 3 & 4 & 1 & 2 & 5 \\
\end{tabular}
\end{center}
\caption{The canonical quandle matrix for $Q_{53}$}
\label{table:CanonicalQ53}
\end{minipage} 
\end{table}

 The next corollary immediately follows 
from Corollary 5 in \cite{HN}, Theorem \ref{thm:automorphism} and Corollary \ref{cor:AtMost}.

\begin{corollary}
  Let $Q$ be a finite connected quandle of order $n$
with its elements naturally ordered,
and with profile $\{ 1, n-1 \}$ or $\{ 1,1,n-2 \}$.
 Then the set of all the automorphisms of $Q$
coincides with the set of all the natural reorderings of $Q$.
\end{corollary}

 There are many quandles as in the above corollary.
 The connected quandles with profile $\{ 1, n-1 \}$ in Vendramin's list 
are shown in Table \ref{table:profile1n-1}.
 For example, the connected quandles with $19$ elements with profile $\{ 1,18 \}$
are $Q_{19,1}$, $Q_{19,2}$, $Q_{19,9}$, $Q_{19,12}$, $Q_{19,13}$ and $Q_{19,14}$.
 The quandle $Q_{62}$ in Vendramin's list has profile $\{ 1,1,4 \}$.

\begin{table}
\begin{center}
\begin{tabular}{|c|c|c|c|c|c|c|c|c|c|c|c|c|c|c|c|}
order    & 1 & 2 & 3 & 4 & 5   & 6 & 7   & 8   & 9   & 10 & 11  & 12 & 13       & 14 & 15 \\
\hline
quandles & 1 &   & 1 & 1 & 2,3 &   & 4,5 & 2,3 & 7,8 &    & 6-9 &    & 1,5,6,10 &    &    \\
\end{tabular}

\begin{tabular}{|c|c|c|c|c|c|c|c|c|c|c|c|c|c|c|}
16  & 17            & 18 & 19          & 20 & 21 & 22 & 23                      & 24 \\
\hline
8,9 & 2,4-6,9-11,13 &    & 1,2,9,12-14 &    &    &    & 4,6,9,10,13,14,16,18-20 &    \\
\end{tabular}

\vspace{1mm}

\begin{tabular}{|c|c|c|c|c|c|c|c|c|c|c|c|c|c|c|}
25    & 26 & 27    & 28 & 29                              & 30 & 31                  \\
\hline
31-34 &    & 62-65 &    & 1,2,7,9,10,13,14,17,18,20,25,26 &    & 2,10-12,16,20,21,23 \\
\end{tabular}

\vspace{1mm}

\begin{tabular}{|c|c|c|c|c|c|c|c|c|c|c|c|c|c|c|}
32     & 33 & 34 & 35 \\
\hline
10-15  &    &    &    \\
\end{tabular}
\end{center}
\caption{\rm Connected quandles with profile $\{ 1, n-1 \}$}
\label{table:profile1n-1}
\end{table}

 Let $Q$ be a finite quandle.
 For an element $i \in Q$,
the map $l_i : Q \ni x \mapsto i*x \in Q$ is not necessarily a bijection in general.
 In \cite{HN},
Ho and Nelson defined a {\it latin} quandle 
to be a quandle with the map $l_i$ being bijection for any $i \in Q$.
 They showed that any conjugation quandle of a group is latin.
 The quandles $Q_{31}$, $Q_{41}$, $Q_{51}$, $Q_{52}$ and $Q_{53}$ in Vendramin's list
are latin.
 $Q_{61}$ shown in Table \ref{table:q61} 
is the first example of a connected non-latin quandle.
 The next theorem is a generalization of Corollary 6.4 in \cite{LR},
and is due to Kanako Oshiro.
 The converse is not true since $Q_{51}$ in Vendramin's list is latin
and of profile $\{ 1,2,2 \}$.

\begin{theorem}\label{thm:latin} {\rm (K. Oshiro)}
 Let $Q$ be a connected quandle with $n$ elements.
 If the profile of $Q$ is $\{ 1,n-1 \}$, then $Q$ is latin.
\end{theorem}

 We prove 
Theorems \ref{thm:automorphism} and \ref{thm:canonical} in the next section,
and Theorem \ref{thm:latin} in Section \ref{sect:latin}.

\section{Proofs of Theorems \ref{thm:automorphism} and \ref{thm:canonical}}
\label{sect:proof}

 Throughout this section,
let $Q$ be a finite connected quandle with $n$ elements $1,2, \cdots, n$
and with profile $\{ 1, \ell_1, \ell_2, \cdots, \ell_k \}$, 
where $1 \le \ell_1 \le \ell_2 \le \cdots \le \ell_k$.

\begin{lemma}\label{lem:ConditionToBeNatural}
 A bijection $\nu : \{ 1,2, \cdots, n\} \rightarrow \{ 1,2,\cdots, n \}$ 
is a natural reordering with respect to $r_{\nu^{-1}(n)}$
if and only if
the new quandle matrix after reordering by $\nu$
is canonical.
\end{lemma}

\begin{proof}
 The $\lq\lq$only if part" is very clear.
 We show the $\lq\lq$if part".
 Let  $r_b$ and $R_b : Q \ni x \mapsto x*b \in Q$
be the permutations 
before and after the reordering by $\nu$ respectively.
 Since the new quandle matrix is canonical,
the permutation $R_n$ is decomposed in the form (N) shown in Introduction.
 Set $\nu^{-1}(n) = q$,
and $\nu^{-1}((\sum_{j=1}^{s-1} \ell_j)+t) = i_{st}$
for any pair of integers $s$ and $t$ with $1 \le s \le k$ and $1 \le t \le \ell_s$.
 Then 
$r_q (i_{st}) 
= i_{st} * q 
= \nu^{-1}((\sum_{j=1}^{s-1} \ell_j)+t)*\nu^{-1}(n)
= \nu^{-1}(((\sum_{j=1}^{s-1} \ell_j)+t)*n)
= \nu^{-1}(R_n(((\sum_{j=1}^{s-1} \ell_j)+t)))
= \nu^{-1}((\sum_{j=1}^{s-1} \ell_j)+(t+1))
= i_{s\,t+1}$
where $t+1$ is read to be an integer in the interval $[1,\ell_s]$ modulo $\ell_s$.
 We can see that the third equality holds
from the way of construction of the new quandle matrix after reordering
shown in Introduction (see Definition 1 in \cite{HN}).
 Hence, 
$\nu$ is the natural reordering
with respect to the presentation
$r_q = (i_{11}\ i_{12}\ \cdots\ i_{1\ell_1})
\cdots ( i_{k1}\ i_{k2}\ \cdots \ i_{k\ell_k} )
(q)$.
\end{proof}

\begin{lemma}\label{lem:q_to_n}
For any element $q$ in $Q$,
there is an automorphism $\nu$ of $Q$
such that $\nu(q) = n$.
\end{lemma}

\begin{proof}
 Since $Q$ is connected,
there is a set of elements $i_1, i_2, \cdots, i_m$ of $Q$
such that 
$((((q \star_1 i_1) \star_2 i_2) \star_3 \cdots) \star_m i_m)=n$,
where $\star_j \in \{ *, \bar{*} \}$ for any integer $j$ with $1 \le j \le m$.
 Then the composition
$\nu=
r_{i_m}^{\epsilon_m} \circ r_{i_{m-1}}^{\epsilon_{m-1}} \circ \cdots \circ r_{i_2}^{\epsilon_2} \circ r_{i_1}^{\epsilon_1}$
with $\epsilon_j=+1$ (when $\star_j=*$) and $-1$ (when $\star_j=\bar{*}$)
brings $q$ to $n$. 
 Note that $\nu$ is an automorphism of $Q$
because $r_{i_j}^{\pm 1}$ is an automorphism of $Q$ for each $j$ with $1 \le j \le m$.
\end{proof}

 In the rest of this section,
we assume that the elements of $Q$ are naturally ordered,
and hence the quandle matrix is in canonical form.
 Under this condition,
we can easily show the next three lemmas.
 We omit the proofs.

\begin{lemma}\label{lem:composite}
 Let $\mu$ be a natural reordering with respect to $r_q$ for some $q \in Q$,
and $\nu$ a natural reordering with respect to $r_n$.
 Then the composition $\nu \circ \mu$ is a natural reordering with respect to $r_q$.
\end{lemma}

\begin{lemma}\label{lem:inverse}
 Let $\nu$ be a natural reordering with respect to $r_n$.
 Then the inverse map $\nu^{-1}$ is also a natural reordering with respect to $r_n$.
\end{lemma}

\begin{lemma}\label{lem:group}
 The set of all the natural reoderings with respect to $r_n$ forms a subgroup
of the symmetric group on $\{ 1,2,\cdots, n \}$.
\end{lemma}

\begin{proof}[Proof of Theorem \ref{thm:automorphism}]
 If $\mu$ is an automorphism of $Q$,
then it fixes the quandle matrix
by Corollary 5 in \cite{HN}.
 In particular, $\mu$ unchanges $r_n$,
and hence the quandle matrix is canonical also after reordering by $\mu$.
 Thus $\mu$ is a natural reordering by Lemma \ref{lem:ConditionToBeNatural}.
\end{proof}

\begin{proof}[Proof of Theorem \ref{thm:canonical}]
 Because $\mu$ is natural with respect to $r_q$, 
we have $\mu(q)=n$.
 Since $Q$ is connected,
there is an automorphism $\lambda$ of $Q$ with $\lambda(q)=n$ 
by Lemma \ref{lem:q_to_n}.
 Then the reordering $\nu= \mu \circ \lambda^{-1}$ fixes $n$.
 Since $\lambda$ is an automorphism of $Q$, 
it fixes the quandle matrix (Corollary 5 in \cite{HN}).
 Hence the new quandle matrices after reordering by $\nu$ and $\mu$ coincide.
 Because the natural reordering $\mu$ fixes $r_n$ in the form (N) in Introduction,
also $\nu$ does,
and hence $\nu$ is a natural reordering with respect to $r_n$ by Lemma \ref{lem:ConditionToBeNatural}.

 Moreover, 
since $\nu$ is natural with respect to $r_n$,
$\nu^{-1}( (\sum_{j=1}^{k-1} \ell_j)+i) =  (\sum_{j=1}^{m-1} \ell_j)+t+i$ holds
for some integers $m, t$ with $1 \le m \le k$, $0 \le t \le \ell_m -1$ and $\ell_m = \ell_k$
and for any integer $i$ with $1 \le i \le \ell_m$,
where $t+i$ is read to be an integer in the interval $[1,\ell_m]$ modulo $\ell_m$.
 Then we consider the reordering 
$\nu'=r_n^t \circ \nu = r_n^t \circ (\mu \circ \lambda^{-1})$.
 Note that
$\nu'(n)=n$ 
and 
$\nu'( (\sum_{j=1}^{m-1} \ell_j)+i)
= ((\sum_{j=1}^{k-1} \ell_j)-t+i) + t 
= (\sum_{j=1}^{k-1} \ell_j)+i$
for any integer $i$ with $1 \le i \le \ell_m$.
 The isomorphism $r_n^t$ does not change the quandle matrix.
 Hence the quandle matrices after reordering by $\mu$ and $\nu'$ coincide,
and $\nu'$ is natural with respect to $r_n$ by Lemma \ref{lem:ConditionToBeNatural}.
\end{proof}

\section{Connected quandle with profile $\{ 1, n-1 \}$}
\label{sect:latin}

 Let $Q$ be a finite connected quandle of order $n$ and with profile $\{ 1, n-1 \}$.

\begin{proof}[Proof of Theorem \ref{thm:latin}]
 Suppose, for a contradiction, that $Q$ is not latin.
 Then $k*i = k*j$ for some elements $i,j,k \in Q$ with $i \ne j$.

 Suppose first that $k=i$ or $k=j$, say $k=i$.
 Then $i*j = k*j = k*i = i*i=i$.
 Hence $r_j$ fixes $j$ and $i$.
 Since $Q$ is of profile $\{ 1, n-1 \}$,
we have $n-1 = 1$, 
and hence $n=2$,
which contradicts the fact that there is no connected quandle of order $2$.

 Then we can assume that $k \ne i$ and $k \ne j$.
 Note that $r_k$ is a cycle of length $n-1$ and fixes $k$,
and hence $r_k$ is a cycle on $n-1$ letters containing $i$ and $j$. 
 Then, there is an integer $m$ with $1 \le m \le n-2$ and $r_k^m (i)=j$.
 Hence,
by Formula (6) in Introduction,
$k*i = k*j = k * (r_k^m (i)) 
= r_k^m ( r_i ( r_k^{-m} (k)))
= r_k^m ( r_i (k))
= r_k^m ( k*i )$.
 This means that $r_k^m$ fixes $k*i $.
 Because $1 \le m \le n-2$,
the cycle $r_k$ fixes $k*i$,
and hence $k=k*i$.
 Then $r_i$ fixes $i$ and $k$,
which is a contradiction.
\end{proof}

\section*{Acknowledgement}
The author would like to thank Kanako Oshiro for helpful comments.


\bibliographystyle{amsplain}

\medskip

\noindent
Department of Mathematical and Physical Sciences,\\
Faculty of Science, Japan Women's University,\\
2-8-1 Mejirodai, Bunkyo-ku, Tokyo, 112-8681, Japan.\\
hayashic@fc.jwu.ac.jp

\end{document}